\documentclass[10pt,twoside]{article}
\usepackage{amssymb}

\newtheorem{thm}{Theorem}[section]

\newtheorem{propn}[thm]{Proposition}

\newtheorem{lmm}[thm]{Lemma}

\newtheorem{cor}[thm]{Corollary}

\newcommand{\R}{{\bf R}}
\newcommand{\Z}{{\bf Z}}

\newcommand{\Q}{{\bf Q}}

\newenvironment{Pf}{\medskip \noindent {\bf Proof: }}
   {$\diamondsuit$ }

\begin{document}
\pagenumbering{arabic}

\begin{center}
\Large{Isometric actions of Heisenberg groups on compact Lorentz manifolds}

\normalsize{\textsc{Karin Melnick}}\footnote{\small{\emph{Department of Mathematics, University of Chicago, Chicago, IL}}}

\end{center}

\begin{quotation}
\small{We prove results toward classifying compact Lorentz manifolds on which Heisenberg groups act isometrically. We give a general construction, leading to a new example, of codimension-one actions---those for which the dimension of the Heisenberg group is one less than the
dimension of the manifold.  The main result is a classification of codimension-one actions, under the assumption they are real-analytic.}
\end{quotation}

\section{Introduction}

Connected isometry groups of compact connected Lorentz manifolds have been classified by Adams and Stuck and independently by Zeghib.  

\begin{thm}[Adams and Stuck \cite{AS}, \cite{AS2}; Zeghib \cite{Ze1}, \cite{Ze2}]
The identity component of the isometry group of a compact connected Lorentz manifold $M$ is isomorphic to $K \times {\R}^k \times S$, where 
\begin{itemize}
\item{$K$ is compact}
\item{$S$, if nontrivial, is one of three types:}
\begin{enumerate}
\item{a finite cover of $PSL_2({\R})$}
\item{a group locally isomorphic to a Heisenberg group $H_n$}
\item{a quotient $S_{\lambda}/ \Lambda$, where $S_{\lambda}$, for each $\lambda \in {\Q}^n$, is isomorphic to $S^1 \ltimes H_n$, and $\Lambda$ is a lattice in the center of $H_n$}
\end{enumerate}
\end{itemize}
\end{thm}

After this classification of isometry groups, the question arises, which manifolds admit these isometry groups.

The groups $PSL_2 ({\R})$ and $S_{\lambda}$ admit bi-invariant Lorentz metrics, so their quotients by cocompact lattices are compact homogeneous Lorentz manifolds.  In \cite{Ze1}, Zeghib shows that these are essentially the only homogeneous compact Lorentz manifolds with noncompact isometry group.  

For an arbitrary Lorentz manifold $M$, if type $1$ in Theorem 1.1 occurs, then Gromov's splitting theorem (\cite{Gromov}, 5.3.D, 5.4.A) says that the universal cover of $M$ is isometric to a warped product $L \times_f S$, where $L$ is a Riemannian manifold, $f$ is a smooth function $L \rightarrow {\R}^+$, and $S$ carries the bi-invariant Lorentz metric coming from the Killing form.  For type $3$, Zeghib has proved a similar splitting theorem (\cite{Ze1} 1.7). We consider here the question of which compact manifolds have isometric actions of groups of type 2.  In this case, the restriction of the metric to the orbits is degenerate, and the orbits cannot be the fibers of a warped product.

Zeghib \cite{Ze2} and Adams have constructed deformations of the $S_{\lambda}$-invariant metric on compact quotients $S_{\lambda}/ \Gamma$ for which some Heisenberg group $H_k$ acts isometrically, but no non-nilpotent group does.  This paper begins with a proof that Adams' deformation on $S_{\lambda} / \Gamma$ has connected isometry group locally isomorphic to the nilradical $H_n$ of $S_{\lambda}$ (Proposition \ref{isomisH}).  We give a general construction of spaces diffeomorphic to $\R \times H_n$ admitting compact quotients with isometric $H_n$-action (Propositions \ref{ConstrIsIsometric} and \ref{GammaIsIsometric}).  There exist such compact quotients not finitely covered by any $S_{\lambda} / \Gamma$.

Our purpose in this paper is to address classification of isometric Heisenberg actions on compact Lorentz manifolds, with particular focus on codimension-one actions---those for which the dimension of the Heisenberg group is one less than the dimension of the manifold.  Our main results are in Section 4.  We classify codimension-one actions, under the assumption they are real-analytic.

{\bf Statement of Results.  } Recall that the $(2n+1)$-dimensional \emph{Heisenberg group}, denoted $H_n$, is the simply-connected Lie group with Lie algebra 
$$\mathfrak{h}_n = < X_1, \ldots, X_n, Y_1, \ldots, Y_n, Z \  | \  [X_i,Y_i] = Z, i=1, \ldots, n>$$

Given a vector space splitting $\mathfrak{h}_n = \mathfrak{p} \oplus \R Z_0$, where $Z_0$ is a generator of the center, there is a non-degenerate symplectic form $\omega_0$ on $\mathfrak{p}$ defined by
$$ [X,Y] = \omega_0(X,Y) Z_0$$

We first describe a necessary condition to have a codimension-one Heisenberg action.  If $H_n$ acts isometrically on a Lorentz manifold $M$ of dimension $2n+2$, then the above splitting determines a vector field $W$ on $M$ transverse to $H_n$-orbits (by conditions ({\bf M1}) - ({\bf M3}) in \S 3.3.1).  In order for $H_n$ to act isometrically, it is necessary that, for any Killing fields $X$ and $Y$ in $\mathfrak{p}$ and $Z \in \R Z_0$,
$$<[W,X],Y>_x  =  \omega_0(X,Y)$$
and
$$ [W,Z](x) = {\bf 0}$$ 
for all $x \in M$. The symmetry of the metric implies 
$$\omega_0(X,[W,Y]) + \omega_0([W,X],Y) = 0$$

for all $X,Y \in \mathfrak{p}$.
It follows that the bracket with $W$ is a derivation of $\mathfrak{h}$.  It must also be \emph{infinitesimally definite} (see \S 3.3), a condition corresponding to positive-definiteness of the metric.
The flow along $W$ from a point $x_0$ of $M$ induces a path $\varphi_t$ in $\mbox{Aut}(\mathfrak{h})$. 

Conversely, given such a path $\varphi_t$ in $\mbox{Aut}(\mathfrak{h})$ for which the path
$$ \varphi_t^{-1} \circ \left( \frac{\partial \varphi_s}{\partial s} \right)_t$$
of derivatives annihilates $Z_0$ and is \emph{infinitesimally definite}, the construction in Section 3.3 gives an $H_n$-invariant Lorentz metric on an $H_n$-space diffeomorphic to $\R \times H_n$.  Provided this path satisfies a compatibility condition with a lattice in $H_n$, the metric descends to a compact quotient $M$ with an isometric $H_n$-action.

We have the following classification (see \S 3.3 for relevant definitions).

{\bf Theorem \ref{classfication}} \emph{There is a bijective correspondence between} 

\emph{(a)  universal covers of compact Lorentz manifolds with real-analytic codimension-one isometric $H_n$-actions, up to equivariant isometry}

\emph{(b) real-analytic, lattice-compatible, metric-defining paths $\Phi_t$ in $\mbox{Aut}(H_n)$, up to equivalence of metric-defining paths.}

The starting point for showing that all codimension-one actions are equivalent to one arising from our construction is Gromov's centralizer theorem (\cite{Gromov} 5.2.A2).  Given a compact real-analytic Lorentz manifold $M$ with isometric $H_n$-action, this theorem gives another copy of $\mathfrak{h}_n$ acting infinitesimally isometrically on the universal cover $\widetilde{M}$, commuting with the lifted $\mathfrak{h}_n$-action.  One of our intermediate results is the following theorem, which describes the possible diffeomorphism types of compact $(2n+2)$-dimensional manifolds $M$ with real-analytic isometric $H_n$-action.

{\bf Theorem \ref{difftype}} \emph{ There is a diffeomorphism of the universal cover $\widetilde{M}$ with $\R \times H_n$ carrying $H_n$-orbits to $H_n$-factors.  There is a lattice $L < H_n$ and a homomorphism $\Psi : \Z \rightarrow \mbox{Aut}(H)$ such that $\pi_1(M) \cong \Z \ltimes_{\Psi} L$. The lattice $L$ acts along the $H_n$-factors of $\widetilde{M}$, and $\Z$ acts by translations along $\R$.}

\bigskip
\emph{Acknowledgements:} I thank my advisor Benson Farb; Abdelghani Zeghib for telling me a key fact for the proof of Proposition \ref{isomisH}; and Scot Adams for discussing with me his construction of Heisenberg actions.  I thank the referee for several helpful comments.

\section{Definitions}

Let $H_n$ be the $(2n+1)$-dimensional Heisenberg group defined in the introduction.  A \emph{warped Heisenberg group} is a solvable extension of $H_n$ with Lie algebra
$$\mathfrak{s}_n = <W, \mathfrak{h} \ | \ [W,X_i] = Y_i, [W,Y_i] = -X_i, i=1, \ldots, n>$$

Note that $\mbox{Ad}(e^{tW})$, with respect to the basis $\{ Z, X_1, Y_1, \ldots, X_n, Y_n, W \}$, is 

$$
\left[
\begin{array}{clrclrc}
1  &   &          &    &    &        &  \\
   & \cos t & - \sin t &    &    &        &   \\
   & \sin t &  \cos t &    &    &        &   \\
   &    &         & \ddots & &       &   \\
   &    &         &    & \cos t & - \sin t&   \\
   &    &         &    & \sin t & \cos t  &   \\
   &    &         &    &        &    & 1
\end{array}
\right]
$$

The group $S_n$ is the group isomorphic to $S^1 \ltimes H_n$ with this Lie algebra of right-invariant vector fields, where $h \mapsto e^{-it}he^{it}$ is the automorphism of $H_n$ with derivative $\mbox{Ad}(e^{tW})$.  

In fact, for any ${\bf \lambda} \in {\R}^n$, there is a warped Heisenberg group with Lie algebra
$$\mathfrak{s}_{\lambda} = <W, \mathfrak{h} \ | \  [X_i,Y_i] = {\lambda}_i Z, [W,X_i] = {\lambda}_i Y_i, [W,Y_i] = - {\lambda}_i X_i, i=1, \ldots, n>$$
 Two such Lie algebras $\mathfrak{s}_{\lambda}$ and $\mathfrak{s}_{\mu}$ are isomorphic if and only if the sets $\{ \lambda_1, \ldots, \lambda_n \}$ and $\{ a \mu_1, \ldots, a \mu_n \}$ are equal for some $a \in \R$ (\cite{AS} 3.2).

Let ${\bf \lambda} \in {\Q}^n$.  In this case, there is a group $S_{\lambda}$ isomorphic to $S^1 \ltimes H_n$ with the Lie algebra $\mathfrak{s}_{\lambda}$.  Here $h \mapsto e^{-it}he^{it}$ has derivative $\mbox{Ad}(e^{tNW})$, where $N$ is the least common multiple of the denominators of the $\lambda_i$. Any cocompact lattice in $H_{\lambda}$ is also cocompact in $S_{\lambda}$. A cocompact lattice $\Gamma$ in $S_{\lambda}$ intersects $H_{\lambda}$ in a uniform lattice and projects modulo $H_{\lambda}$ to a finite subgroup of $S^1$ (see \cite{Rag} 3.3).  Also, a cocompact lattice in $H_{\lambda}$ intersects the center $Z(H_{\lambda}) = Z(S_{\lambda})$ in a cocompact subgroup (see \cite{Rag} 2.3).  

For any ${\bf \lambda}$, the Lie algebra $\mathfrak{s}_{\lambda}$ admits a bi-invariant inner product of signature $(1,2n+1)$ in which the vectors $X_1, \ldots, X_n, Y_1, \ldots, Y_n$ are orthonormal; the vectors $Z$ and $W$ are both isotropic and orthogonal to $X_i$ and $Y_i$ for $i = 1, \ldots, n$; and $<W,Z> = 1$.  This inner product gives rise to a bi-invariant Lorentz metric on $S_{\lambda}$.  If $\Gamma$ is a lattice in $S_{\lambda}$, the metric descends to the quotient $S_{\lambda} / \Gamma$, for which the isometry group has identity component $S_{\lambda}/ \Lambda$, where $\Lambda = \Gamma \cap Z(S_{\lambda})$.  This last fact can be proved by the methods of Propositon \ref{isomisH} below.  Note that the center $Z(S_{\lambda})$ acts here via a compact quotient isomorphic to $S^1$.

\section{Construction of $H$-actions}

Section 3.1 contains a brief account of Zeghib's method for deforming the homogeneous metric on compact manifolds $S_n / \Gamma$, and Section 3.2 treats Adams' deformation of such metrics.  Proposition \ref{isomisH} says that for Adams' deformation, the resulting connected isometry group is $H_n/\Lambda$, where $\Lambda$ is a lattice in $Z(H_n)$.  In Section 3.3 a more general construction of codimension-one Heisenberg actions, leading to a new example, is given.
  
\subsection{Perturbing the holonomy for $S / \Gamma$}

In (\cite{Ze2} 4.1.2) Zeghib constructs metrics on compact quotients of $S_n$ with connected isometry groups containing ${\R}^k$ or $H_k$, $k<n$.  Given $\Gamma$ a lattice in $S_n$ and a homomorphism $\rho: \Gamma \rightarrow S_n$, the group $\Gamma_{\rho} = \{ (\rho(\gamma), \gamma) : \gamma \in \Gamma \}$ is a discrete subgroup of Isom$(S)$.  It is possible to choose $\rho$ so that $\Gamma_{\rho}$ still acts freely and properly on $S$.  Then the connected component of the isometry group of $S / \Gamma_{\rho}$ is the centralizer of $\Gamma_{\rho}$.  One possibility for the Zariski closure Zar($\rho(\Gamma))$ is an $(n-k)$-dimensional abelian subgroup of $H$, where $k < n$, having centralizer in $S_n$ isomorphic to ${\R}^{n-k} \times H_k$.  If the centralizer of $\Gamma_{\rho}$ is non-abelian, then $\rho(Z(\Gamma)) = 1$.  In this case, $Z(S_n) \cap \Gamma_{\rho} = Z(\Gamma)$, so the center of the isometry group of $S_n / \Gamma_{\rho}$ is one-dimensional and acts via a compact quotient.

\subsection{Breaking the symmetry on $S/\Gamma$}

Assume for simplicity that $\lambda = (1, \ldots, 1)$, and denote $H_n$ and $S_n$ simply by $H$ and $S$, respectively.  Also assume $\Gamma \subseteq H$.

Let $m$ be the pullback of a positive bump function at $1$ on $S^1  = H \backslash S$ to $S$.  This function is left-$H$-invariant, right-$\Gamma$-invariant, and not invariant under left translation by any $s \notin H$.  Let $\theta$ be the bi-invariant metric on $S$, and denote by $m \theta$ this metric rescaled by $m$.  Then $m \theta$ descends to a metric on $S/\Gamma$, and $H$ acts here by isometries.  
Represent elements of the semi-direct product $S$ by ordered pairs $(t,h), 0 \leq t < 2 \pi$.  

Since Killing fields on $S/ \Gamma$ coming from the $H$-action do not have constant norm in $m \theta$, it is not immediate that the orbits of one-parameter subgroups in $H$ are geodesics.  The fact that they are geodesic can be seen using results of Zeghib's [Ze3] paper on totally geodesic lightlike foliations in Lorentz manifolds.  

A \emph{lightlike subspace} of a Lorentz inner-product space is one for which the restriction of the inner product is degenerate.  In this case, the kernel is one-dimensional, and the inner product restricted to any complementary subspace is definite.  A \emph{lightlike submanifold} $N$ of a Lorentz manifold $M$ is one for which $T_x N$ is a lightlike subspace of $T_x M$ for all $x \in N$.  In Section 1.1 of \cite{Ze3}, Zeghib shows that a codimension-one lightlike submanifold $N$ of a compact Lorentz manifold is totally geodesic if and only if the kernel distribution on $N$ is \emph{transversally Riemannian}---the flow generated by a vector field tangent to this distribution preserves the degenerate Riemannian metric on $N$.

The following lemma is key in classifying codimension-one actions.

\begin{lmm}  
\label{lemma}
A linear transformation that preserves a Lorentz inner product and restricts to the identity on a codimension-one lightlike subspace is the identity. 
\end{lmm}

\begin{Pf}
Suppose $V$ is a Lorentz inner-product space and $f \in O(V)$ is the identity on $U \subset V$ as above.  Let $Z, X_1, \ldots, X_n$ be a basis for $U$ such that $Z$ spans the kernel of $U$ and $X_1, \ldots, X_n$ are orthonormal.  The orthogonal complement to  $span \{ X_1, \ldots, X_n \}$ is Lorentz and contains a unique isotropic vector $W$ with $<W,Z> = 1$.  Clearly, if $f$ fixes $U$ pointwise, then it must fix $W$.
\end{Pf}

\begin{propn} $\mbox{Isom}^0 (S / \Gamma, m \theta) = H/\Lambda$, for $\Lambda = Z(H) \cap \Gamma$.
\label{isomisH}
\end{propn}

\begin{Pf}
It is clear from the construction that $H$ acts isometrically on $(S/ \Gamma, m \theta)$, and that $\Lambda = Z(H) \cap \Gamma$ is the kernel of this action.  It suffices to show any isometry in the component of the identity lies in $H / \Lambda$.  Suppose $\overline{f} \in \mbox{Isom}^0 (S/\Gamma, m \theta)$.  Note $\overline{f}_*$ is the identity on $\pi_1 (S / \Gamma)$.  By post-composing with an isometry in $H/ \Lambda$ if necessary, one may assume $\overline{f}$ takes the $\Gamma$-coset of $(0,e)$ to the coset of
$(T,e)$ for some $T$.  Let $f$ be the lift of $\overline{f}$ to $S$ mapping $(0,e)$ to
$(T,e)$.  Note that $f$ is equivariant with respect to right translation by $\Gamma$.  

Each $H$-fiber $H_t = \{ t \} \times H$ is lightlike with kernel distribution
generated by the Killing field corresponding to a generator of the center
of $H$.  Thus this distribution is transversally Riemannian and every $H_t$ is
totally geodesic.  In this case, $H_t$ inherits a connection from $S$.  Let $X^* \in \mathfrak{h}$, and let $X$ be the Killing field on $S$ coming from left translation by $e^{tX^*}$.  The squared norm $<X,X> = m \theta(X,X)$ is constant along $H_t$, so for any vector $Y$ tangent to $H_t$, 
$$Y<X,X> = 0 = 2<\nabla_Y X,X>$$

  Pick a point $p \in H_t$ and extend $Y$ along the curve $e^{tX^*} p$ such that $[X,Y] = 0$. Since the flow along $X$ is an isometry, $Y$ also satisfies $X<X,Y> = 0$.  The connection is torsion-free, so at $p$, 
\begin{eqnarray*}
 <\nabla_X X,Y> + <X, \nabla_X  Y> & =  & <\nabla_X X, Y> + < X,\nabla_Y X> + <X, [X, Y]> \\           & =  & <\nabla_X X, Y> \\
           & =  & 0
\end{eqnarray*}

The tangent vector $Y$ to $H_t$ was arbitrary, so an integral curve for $X$ is geodesic in the induced connection on $H_t$, and hence in $S$.  Thus, all the Killing fields in $S$ coming from left translation by elements of $H$ integrate to geodesics.

Identify the tangent space to an $H$-fiber $T_{(t,e)} H_t$ with $\mathfrak{h}$ by evaluating Killing fields at $(t,e)$.  At each $(t,e)$ there are two restricted exponential maps $ T_{(t,e)} H_t \rightarrow
H_t$, one coming from the exponential map in $H$ and one from the
connection.  These must agree; denote the map by $\mbox{exp}_t$.  Since the Lie group exponential on a simply connected nilpotent Lie group is a diffeomorphism (see \cite{Knapp} 1.104, \cite{Rag} 1.9), $\mbox{exp}_t$ maps $T_{(t,e)} H_t$ diffeomorphically onto the fiber $H_t$.  Because $\Gamma \cap Z(H)$ is a lattice in $Z(H)$, the inverse image $\mbox{exp}_0^{-1}(\Gamma \cap Z(H))$ spans the center $\mathfrak{z} \subset T_{(0,e)} H_0$.  Then any subspace of $\mathfrak{h}$ containing $\mbox{exp}_0^{-1}(\Gamma)$ is a subalgebra.  The corresponding closed connected subgroup of $H$ must be cocompact, so it must equal $H$.  Thus $\mbox{exp}_0^{-1}(\Gamma)$ contains a basis of $\mathfrak{h}$, and an isometry is determined along $H_0$ by its derivative on this subset.  Now $\Gamma$-equivariance implies that $f_{*(0,e)}: T_{(0,e)} H_0 \rightarrow T_{(T,e)} H_T$ maps $\mbox{exp}_0^{-1}(0,\gamma)$ to $\mbox{exp}_T^{-1} (T,\gamma)$ for all $\gamma \in \Gamma$.  The element of $\mathfrak{h}$ evaluating to $\mbox{exp}_T^{-1}(T,\gamma)$ at $(T,e)$ is $\mbox{Ad}(e^{-TW})(\mbox{exp}_0^{-1}(\gamma)$.  Thus $f_{*(0,e)}$ induces Ad$(e^{-TW})$ on $\mathfrak{h}$, and $f$ maps $H_0$ to $H_T$ by $(0,h) \mapsto (T,h)$.  But $f$ is an isometry, so $m(0) = m(T)$.  Hence $0 = T$, and $f_{*(0,e)}$ restricts to the identity on $T_{(0,e)} H_0$.  

By the lemma \ref{lemma}, $f_{*(0,e)} \equiv \mbox{Id}$.  An isometry of a connected pseudo-Riemannian manifold fixing a point and a frame is trivial (see \cite{Ko} II.1.3).  
\end{Pf}

\subsection{More general codimension-one actions}

\subsubsection{Construction of $\R *_{\Phi} H$}

Fix a vector space splitting $\mathfrak{h} = \mathfrak{p} \oplus \R Z_0$, where $Z_0 \neq {\bf 0}$ belongs to the center $\mathfrak{z}$.  Let $\omega_0$ be the non-degenerate symplectic form on $\mathfrak{p}$ determined by $Z_0$:
$$[X,Y] = \omega_0(X,Y) Z_0$$

Call an endomorphism $\nu$ of $\mathfrak{p}$ \emph{infinitesimally definite} if $\omega_0(X,\nu(X)) > 0$ for all $X \in \mathfrak{p}$.  Note that an infinitesimally definite endomorphism is invertible with $\omega_0(\nu^{-1}(X),X) > 0$ for all $X$.  Further, $\nu$ is \emph{infinitesimally symplectic}, meaning
$$ \omega_0(\nu(X),Y) + \omega_0(X,\nu(Y)) = 0$$

if and only if $\nu^{-1}$ is.  An endomorphism $\delta$ of $\mathfrak{h}$ annihilating $Z_0$ is infinitesimally symplectic on $\mathfrak{h} / \mathfrak{z}$ if and only if it is a derivation.  Let $\mathcal{M}_{\mathfrak{p}}$ be the set of derivations $\delta$ of $\mathfrak{h}$ preserving $\mathfrak{p}$ such that $\delta(Z_0) = {\bf 0}$ and $\nu = \delta |_{\mathfrak{p}}$ is infinitesimally definite. 
Note that $\mathcal{M}_{\mathfrak{p}}$ is not a subalgebra of $\mbox{Der}(\mathfrak{h})$ and does not depend on the choice of $Z_0$.  Call a path $\varphi : \R \rightarrow \mbox{Aut}(\mathfrak{h})$ \emph{metric-defining} if for some $\mathfrak{p}$ 
$$\varphi_t^{-1} \circ \left( \frac{\partial \varphi_s}{\partial s} \right)_t  \in \mathcal{M}_{\mathfrak{p}} \qquad \mbox{for all} \ t$$  

Call a path $\Phi : \R \rightarrow \mbox{Aut}(H)$ \emph{metric-defining} if the path $\varphi_t = D_e (\Phi_t)$ in $\mbox{Aut}(\mathfrak{h})$ is metric-defining.

Given a path $\Phi_t$ in $\mbox{Aut}(H)$, define an $H$-space $\R *_{\Phi} H$ diffeomorphic to $\R \times H$ by
$$ h.(t,g) = (t, \Phi_t(h) g)$$

This $H$-action is free, so there exist at each $x$ linear injections 
\begin{eqnarray*}
f_x & : & \mathfrak{h} \rightarrow T_x(\R *_{\Phi} H) \\ 
 f_x(Y) & = & \left. \frac{\partial}{\partial t} \right|_0 (e^{tY}.x)
\end{eqnarray*}

Suppose that $\Phi_t$ is metric-defining, and let $\varphi_t = D_e(\Phi_t)$.  There is $\mathfrak{p}$ such that

$$  \varphi_t^{-1} \circ \left( \frac{\partial \varphi_s}{\partial s} \right)_t   \in \mathcal{M}_{\mathfrak{p}} \ \mbox{for all} \ t$$
 
Note that $\mathfrak{p}$ is determined by $\varphi_t$ because it is exactly the image of the derivation above.  Denote by $\nu_t$ the restriction of this derivation to $\mathfrak{p}$.  Let $W$ be the vector field $\frac{\partial}{\partial t}$ on $\R *_{\Phi} H$. Define a Lorentz metric $<,>_x$ on $\R *_{\Phi} H$ in which $W$ satisfies at each $x$:
\begin{eqnarray*}
  & <W(x),W(x)>_x = 0 & ({\bf M1}) \\ 
  & <W(x),f_x(Z_0)>_x = 1 & ({\bf M2}) \\
  & W(x) \perp f_x(\mathfrak{p}) &  ({\bf M3})
\end{eqnarray*}

Further, make $f_x(Z_0)$ isotropic and $f_x(Z_0) \perp f_x(\mathfrak{p})$ everywhere.  Last, for $X,Y \in \mathfrak{p}$, let
$$ <f_x(X),f_x(Y)>_x = \omega_0( \nu_{t(x)}^{-1}(X), Y)$$

where $t(x)$ is the $\R$-coordinate of $x$. Since $\nu_{t(x)}$ is the restriction of a derivation of $\mathfrak{h}$, both $\nu_{t(x)}$ and $\nu_{t(x)}^{-1}$ are infinitesimally symplectic.  It follows that the form $<,>_x$ is symmetric.  Because $\nu_{t(x)}$ is infinitesimally definite, the form is positive definite on $\mathfrak{p}$.  Note that this metric on $\R *_{\Phi} H$ is associated to the path $\Phi_t$ and the generator $Z_0$.  In this section, the notation $(\R *_{\Phi} H, Z_0)$ will stand for the Lorentz manifold thus determined.

\begin{propn}
\label{ConstrIsIsometric}
The metric on $\R *_{\Phi} H$ is $H$-invariant.
\end{propn}

\begin{Pf}
The $H$-action is isometric if and only if, for each $K \in \mathfrak{h}$, $x \in \R *_{\Phi} H$, and vector fields $X,Y$ near $x$,
\begin{eqnarray*}
0 & = & \left. \frac{\partial}{\partial t} \right|_0 < e^{tK}_* X,e^{tK}_* Y>_{e^{tK} x} \\
 & = & <[f_x(K),X],Y>_x + <X,[f_x(K),Y]>_x + f_x(K)<X,Y>
\end{eqnarray*}

Since the above sum depends on $X$ and $Y$ only for their values at $x$, it suffices to check vanishing for $X(x)$ and $Y(x)$ ranging over a basis of the tangent space at $x$.  For $X \in \mathfrak{h}$, denote the vector field given by $f_x(X)$ simply by $X^*$.  Since $X$ corresponds to a right-invariant vector field on $H$, and $H$ acts on the left, $X \mapsto X^*$ is a Lie algebra homomorphism.  Take $Y_1, \ldots, Y_{2n}$ a basis for $\mathfrak{p}$.  If $X$ and $Y$ each coincide near $x$ with elements of $\{ Z_0^*, Y_1^*, \ldots, Y_{2n}^*, W \}$, then $K^*(x) <X,Y> = 0$ for any $K \in \mathfrak{h}$, because any such inner product is constant along $H$-orbits.  Now the equation is
$$<[K^*,X],Y>_x + <X,[K^*,Y]>_x = 0$$

for all $K \in \mathfrak{h}$ and $X,Y \in \{ Z_0^*, Y_1^*, \ldots, Y_{2n}^*, W \}$.  If $X$ and $Y$ both come from $\mathfrak{h}$, then both terms above vanish because $Z_0^*(x) \perp f_x(\mathfrak{h})$.

It remains to check vanishing when $X = W$.  Let the coordinates of $x$ be $t(x) = t$ and $p(x) = p$.  
The flow along $K^*$ is 
\begin{eqnarray*}
e^{rK}.x & = & (t,\Phi_t(e^{rK}) \cdot p) \\
& = & (t,e^{r \varphi_t(K)} \cdot p)
\end{eqnarray*}

The derivative on $W(x)$ is 
$$ e^{rK}_*(W(x)) = (\frac{\partial}{\partial t}, r \left( \frac{\partial \varphi_s}{\partial s} \right)_t (K) \cdot (e^{r \varphi_t(K)} \cdot p))$$

where $X \cdot h$, for $X \in \mathfrak{h}$ and $h \in H$, denotes the image of $X$ under the derivative of right translation by $h$.  Finally, the bracket 
\begin{eqnarray*}
[K^*,W](x) & = & \left. \frac{\partial}{\partial r} \right|_0 e^{-rK}_*(W(e^{rK}.x)) \\
& = & \lim_{r \rightarrow 0} \frac{1}{r} (e^{-rK}_*(W(e^{rK}.x)) - W(x))  \\
& = & ( {\bf 0}, \lim_{r \rightarrow 0} \frac{1}{r} (- r \left( \frac{\partial \varphi_s}{\partial s} \right)_t (K) \cdot p) \\
& = & ({\bf 0}, - \left( \frac{\partial \varphi_s}{\partial s} \right)_t (K) \cdot p) \\
& = & - (\varphi_t^{-1} \circ \left( \frac{\partial \varphi_s}{\partial s} \right)_t (K))^*(x)
\end{eqnarray*}

If $K = Z_0$, then 
$$[Z_0^*,W] = -(\varphi_t^{-1} \circ \left( \frac{\partial \varphi_s}{\partial s} \right)_t)  (Z_0)^* = {\bf 0}$$

If $K \in \mathfrak{p}$, then we have
$$<[K^*,W],Y>_x + <W,[K^*,Y]>_x = <- \nu_t(K)^* ,Y>_x + <W,[K^*,Y]>_x$$

When $Y = Z_0^*$, then both terms vanish because $\nu_t(K)^* \in f_x(\mathfrak{p})$, which is orthogonal to $Z_0^*(x)$, and $[K^*,Z_0^*] = {\bf 0}$.  When $Y= Y^*$ for $Y \in \mathfrak{p}$, then we have
\begin{eqnarray*}
<- \nu_t(K)^*,Y^*>_x + <W,[K^*,Y^*]>_x & = & <- \nu_t(K)^*,Y^*>_x + <W,\omega_0(K,Y) Z_0^*>_x \\
 & = & - <\nu_t(K)^*,Y*>_x + \omega_0(K,Y) \\
& = & 0
\end{eqnarray*}

since $<X^*,Y^*>_x = \omega_0(\nu_t^{-1} X,Y)$ for all $X \in \mathfrak{p}$.
\end{Pf}

\subsubsection{Equivalence}

Next we consider when the spaces associated to pairs $(\Phi_t,Z_0)$ and $(\Phi_t^{\prime},Z_0^{\prime})$ are $H$-equivariantly isometric.

Suppose given a generator $Z_0 \in \mathfrak{z}$ and a path $\Phi_t$, metric-defining with respect to a complementary subspace $\mathfrak{p}$.  Let $\varphi_t = D_e(\Phi_t)$ and $\nu_t$ be as usual.  Suppose $\mathfrak{p}^{\prime}$ is a different complement to $\mathfrak{z}$.  There is a canonical isomorphism $a : \mathfrak{p} \rightarrow \mathfrak{p}^{\prime}$ inducing the identity on $\mathfrak{h} / \mathfrak{z}$.  This isomorphism has the form 
$$ a(X) = X + \alpha(X) Z_0$$ 

for some functional $\alpha$ on $\mathfrak{p}$.  Let $A \in \mathfrak{p}$ be such that $\omega_0(A,X) = \alpha(X)$ for all $X \in \mathfrak{p}$.

  Define
$$ \nu_t^{\prime} = a \circ \nu_t \circ a^{-1}$$

and let $\varphi_t^{\prime}$ and $\Phi_t^{\prime}$ be corresponding metric-defining paths.  Let $\gamma$ be a path in $\mathfrak{h}$ satisfying
$$\dot{\gamma}(t) = \left( \frac{\partial \varphi_s^{\prime}}{\partial s} \right)_t (A) - \frac{1}{2} \alpha( \nu_t^{\prime} (A) ) Z_0 $$

Then 
$$ \theta : (t,h) \mapsto (t, \Phi_t^{\prime} \circ \Phi_t^{-1}(h) \cdot e^{\gamma(t)})$$

is an $H$-equivariant isometry between $(\R *_{\Phi} H, Z_0)$ and $(\R *_{{\Phi}^{\prime}} H, Z_0)$.  

The $H$-equivariance is clear.  By equivariance, it suffices to check isometry at $(t,e)$ for all $t$; also, Killing fields are carried to Killing fields.  Denote by $<,>$ and $<,>^{\prime}$ the metrics on the domain and target, respectively.  The restriction of $\theta_{*x}$ to $f_x(\mathfrak{h})$ is isometric: for any $X, Y \in \mathfrak{p}$, 
\begin{eqnarray*}
<X^*,Y^*>_x & = & \omega_0(\nu_{t(x)}^{-1} X,Y) \\
& = & \omega_0((a \circ \nu_{t(x)}^{-1} \circ a^{-1}) (a(X)),a(Y) ) \\
& = & \omega_0((\nu_{t(x)}^{\prime})^{-1}(a(X)),a(Y)) \\
& = & <a(X)^*,a(Y)^*>^{\prime}_{\theta(x)} \\
& = & <X^*,Y^*>^{\prime}_{\theta(x)}
\end{eqnarray*}

Denote by $W$ and $W^{\prime}$ the vector fields $\frac{\partial}{\partial t}$ on the domain and target, respectively.  
\begin{eqnarray*}
 \theta_{*(t,e)} (W) & = & ( \frac{\partial}{\partial t}, \dot{\gamma}(t) \cdot e^{\gamma(t)})   \\
& = & W^{\prime} + (\varphi_t^{\prime})^{-1} (\dot{\gamma}(t) )^* \\
& = & W^{\prime} + \nu_t^{\prime}(A)^* - \frac{1}{2} \alpha(\nu_t^{\prime}(A)) Z_0^*
\end{eqnarray*}
using that $\varphi_t^{\prime}$ fixes $Z_0$.  Now it is easy to check $\theta_{*(t,e)}(W)$ is isotropic, orthogonal to $f_{\theta(0,e)}(\mathfrak{p})$, and has inner product 1 with $Z_0^*$.  

Finally, note that for arbitrary $d \in \R$, one may set $\nu_{t+d}^{\prime} = a \nu_t a^{-1}$, and then
$$(t,h) \mapsto (t+d, \Phi_{t+d}^{\prime} \circ \Phi_t^{-1}(h) \cdot e^{\gamma(t)})$$

is an $H$-equivariant isometry of $(\R *_{\Phi} H, Z_0)$ with $(\R *_{\Phi^{\prime}} H, Z_0)$, where $\Phi^{\prime}_t$ is a metric-defining path arising from $\nu_t^{\prime}$, and $\gamma(t)$ is defined in an obvious way.

Next suppose that a metric-defining path $\Phi_t$ and a generator $Z_0$ are given, and define $\varphi_t$ and $\nu_t$ as usual.  Given $c \neq 0$ and $d \in \R$, let 
$$c \nu_{ct+d}^{\prime} = \nu_t$$

Let $\varphi_t^{\prime}$ and $\Phi_t^{\prime}$ be metric-defining paths corresponding to $\nu_t^{\prime}$.  Note that $\Phi_{ct+d}^{\prime} = \Phi_t$.  Then
$$ (t,h) \rightarrow (ct+d,\Phi_{ct+d}^{\prime} \circ \Phi_t^{-1}(h)) = (ct+d,h)$$

is an $H$-equivariant isometry between $(\R *_{\Phi} H, Z_0)$ and $(\R *_{\Phi^{\prime}} H, cZ_0)$.  We leave the verification to the reader.

These are the only ways to obtain $H$-equivariantly isometric spaces.

\begin{propn}
\label{HEquiv}
Suppose there is an $H$-equivariant isometry 
$$ \theta : (\R *_{\Phi} H, Z_0) \rightarrow (\R *_{\Phi^{\prime}} H, Z_0^{\prime})$$

Let $\nu_t$ and $\nu_t^{\prime}$ be the restricted derivations of $\mathfrak{h}$ associated to $\Phi_t$ and $\Phi_t^{\prime}$, respectively.  Let $\mathfrak{p}$ and $\mathfrak{p}^{\prime}$ be their respective images, and $a$ the canonical isomorphism $\mathfrak{p} \rightarrow \mathfrak{p}^{\prime}$.
Then there exist $c \neq 0$ and $d \in \R$ such that $Z_0^{\prime} = cZ_0$, and $c \nu_{ct+d}^{\prime} = a \circ \nu_t \circ a^{-1}$.
\end{propn}

\begin{Pf}
Since $\theta$ preserves $H$-orbits, it has the form
$$\theta: (t,h) \mapsto (H(t),K(t,h))$$

By $H$-equivariance, $\theta_*$ carries each Killing field $X^*$ on the domain to $X^*$ on the target.  Let $W$ and $W^{\prime}$ be the vector fields $\frac{\partial}{\partial t}$ on the domain and target, respectively.  Let $t(x) = t$.  Since $<W, Z_0^*> = 1$ and $Z_0^* \perp f_y(\mathfrak{h})$ for all $y$,
$$<\theta_*(W),Z_0^*>_{\theta(x)} \ = \ <\dot{H}(t) W^{\prime}, Z_0^*>_{\theta(x)}\  = \ 1$$

There must exist $c \neq 0$ such that $Z_0^{\prime} = c Z_0$ and $\dot{H}(t) = c$ for all $t$.  Then $H(t) = ct +d$ for some $d$.  

Denote by $\omega_0$ and $\omega_0^{\prime}$ the symplectic forms defined by $Z_0$ and $Z_0^{\prime}$, respectively.  For $X,Y \in \mathfrak{p}$,
$$<[W,X^*],Y^*>_x = \omega_0(\nu_t^{-1} X,Y) $$

which must equal
\begin{eqnarray*}
<[\theta_*(W),X^*],Y^*>_{\theta(x)} & = & <[ \theta_*(W),a(X)^*],a(Y)^*>_{\theta(x)}  \\
& = & \omega_0^{\prime}((\nu_{ct+d}^{\prime})^{-1}(a(X)),a(Y)) \\
& = & \frac{1}{c} \omega_0((a^{-1} \circ (\nu^{\prime}_{ct+d})^{-1} \circ a) (X), Y)
\end{eqnarray*} 

 Since $\omega_0$ is non-degenerate on $\mathfrak{p}$, we have 
$$\frac{1}{c} a^{-1} \circ (\nu^{\prime}_{ct+d})^{-1} \circ a = \nu_t^{-1},$$

so
$$ c \nu_{ct+d}^{\prime} = a \circ \nu_t \circ a^{-1}$$
\end{Pf}

We will say that two metric-defining paths $\Phi_t$ and $\Phi_t^{\prime}$ are \emph{equivalent} if there are $c \neq 0$ and $d \in \R$ such that $c \nu_{ct+d}^{\prime} = a \circ \nu_t \circ a^{-1}$, where $a$ is the canonical isomorphism as above.

\subsubsection{Compact quotients}

Suppose that $\Phi_t$ is \emph{$\Z$-equivariant}: $ \Phi_{q+t} = \Phi_q \Phi_t$ for all $t \in \R, \ q \in \Z$;
in particular, $\Phi$ restricts to a homomorphism on $\Z$.  Say that $\Phi$ is \emph{lattice-compatible} if it is $\Z$-equivariant and there is some lattice $L$ of $H$ preserved by $\Phi(q)$ for all $q \in \Z$.  Then the semi-direct product $\Gamma = \Z \ltimes_{\Phi} L$ acts on $\R *_{\Phi} H$ by
$$ (t,g).(q,\lambda) = (t + q, \Phi_q(g)\cdot \lambda)$$

The action is written on the right to emphasize that it commutes with the $H$-action.  In particular, it preserves the Killing fields $X^*$ for each $X \in \mathfrak{h}$.  Thus $H$ acts on the compact quotient $M = (\R *_{\Phi} H) / \Gamma$.

\begin{lmm}
\label{PeriodicIfEquiv}
If the paths $\Phi_t$ and $\varphi_t$ are $\Z$-equivariant, then $\nu_t$ is $\Z$-invariant.
\end{lmm}

\begin{Pf}
 Let $q \in \Z$. The endomorphism $\nu_{q+t}$ is the restriction to $\mathfrak{p}$ of
\begin{eqnarray*}
\varphi_{q+t}^{-1} \circ \left( \frac{\partial \varphi_s}{\partial s} \right)_{q+t}  & = & \varphi_{q+t}^{-1} \circ \left( \frac{\partial \varphi_{q+s}}{\partial s} \right)_t  \\
& = &  (\varphi_q \circ \varphi_t)^{-1} \circ \left( \frac {\partial (\varphi_q \circ \varphi_s)}{\partial s} \right)_t  \\
& = & \varphi_t^{-1} \circ \varphi_q^{-1} \circ \varphi_q \circ \left( \frac{\partial \varphi_s}{\partial s} \right)_t  \\
& = & \varphi_t^{-1} \circ  \left( \frac{\partial \varphi_s}{\partial s} \right)_t   
\end{eqnarray*}

which restricts to $\nu_t$ on $\mathfrak{p}$.
\end{Pf}

\begin{propn}
\label{GammaIsIsometric}
The $\Gamma$-action is isometric, so the metric descends to $M$.
\end{propn}

\begin{Pf}
The vector fields $W$ and $Z_0^*$ are $\Gamma$-invariant, so $\Gamma$ is isometric restricted to their span at any $x$.  Also, $f_x(\mathfrak{p})$ is a $\Gamma$-invariant distribution, so it suffices to check that for each $\gamma \in \Gamma$ and $X,Y \in \mathfrak{p}$,
$$<\gamma_* (X^*(x)), \gamma_* (Y^*(x))>_{\gamma x} = <X^*,Y^*>_x$$ 

at all $x \in M$.  Denote by $t(x)$ the $\R$-coordinate of $x$. For any $\gamma \in \Gamma$,
\begin{eqnarray*}
<\gamma_* (X^*(x)), \gamma_* (Y^*(x))>_{\gamma x} & = & <X^*(\gamma x),Y^*(\gamma x)>_{\gamma x} \\
     & = & \omega_0( \nu_{t(\gamma x)}^{-1} (X), Y) \\
     & = & \omega_0(\nu_{t(x)}^{-1}(X),Y) \\
     & = & < X^*,Y^*>_x
\end{eqnarray*}

using $\Z$-equivariance of $\Phi_t$ and Lemma \ref{PeriodicIfEquiv}.
\end{Pf}

\subsubsection{Examples}

{\bf Example.}  Adams' deformation.

Let $W, X_1, \ldots, X_n, Y_1, \ldots, Y_n, Z_0$ be the basis for $\mathfrak{s}$ and $\theta$ be the bi-invariant Lorentz metric on $\widetilde{S}$ from \S 2.  Let $\mathfrak{p}$ be the span of $X_1, \ldots, Y_n$. Let $m$ be the positive function on $S^1$ from \S 3.2.  Denote again by $m$ the lift of this function to $\widetilde{S} \cong \R \ltimes H$.  Let 
$$H(t) = \int_0^t m(s) \mathrm{d}s$$

Note that $H$ is invertible.  Let $\varphi_{H(t)}$ be the automorphism of $\mathfrak{h}$ fixing $Z_0$ and acting on each subspace $span \{ X_i, Y_i \}$ by
$$
\left[
\begin{array}{lr}
\cos t & - \sin t \\
\sin t & \cos t
\end{array}
\right]
$$

Note that $\varphi_{H(t)} = \mbox{Ad}(e^{tW})$.  Now $\left( \frac{\partial \varphi_s}{\partial s} \right)_{H(t)}$ annihilates $Z_0$ and acts on each $span \{ X_i, Y_i \}$ by 
$$
m(t)^{-1} \left[
\begin{array}{lr}
- \sin t & - \cos t \\
\cos t & - \sin t
\end{array}
\right]
$$

Let $\nu_{H(t)}$ be the restriction to $\mathfrak{p}$ of $\varphi_{H(t)}^{-1} \circ \left( \frac{\partial \varphi_s}{\partial s} \right)_{H(t)} $; it acts on each $span \{ X_i, Y_i \}$ by
$$
\left[
\begin{array}{lr}
  &  - m(t)^{-1} \\
m(t)^{-1} & 
\end{array}
\right]
$$

These are infinitesimally definite derivations for each $t$, so $\varphi$ is a metric-defining path.  Let $\Phi_{H(t)}$ be the automorphism of $H$ with derivative $\varphi_{H(t)}$.  Let $\rho$ be the $H$-invariant Lorentz metric $\R *_{\Phi} H$ associated with $Z_0$. The map 
\begin{eqnarray*}
F & : & (\widetilde{S},m \theta) \rightarrow (\R *_{\Phi} H, \rho) \\
F & : & (e^{tW},h) \mapsto (H(t), \Phi_{H(t)} (e^{tW} h e^{-tW})) = (H(t),h) 
\end{eqnarray*}

is an $H$-equivariant isometry; we leave the verifications to the reader.  

Let $\alpha = H(2 \pi)$.  Periodicity of $m$ implies $H(2n \pi + t) = n \alpha + H(t)$.  The automorphisms $\Phi_{n\alpha}$ are trivial for each $n \in \Z$.  To see $\Phi$ is ($\alpha \Z$)-equivariant
\begin{eqnarray*}
\Phi_{n \alpha + H(t)} & = & \Phi_{H(2n \pi + t)} \\
& = & \Phi_{H(t)} \\
& = & \Phi_{n \alpha} \circ \Phi_{H(t)}
\end{eqnarray*}

For any lattice $L$ in $H$, the group $(\alpha \Z) \ltimes_{\Phi} L$ acts isometrically on $\R *_{\Phi} H$, commuting with the $H$-action. 

{\bf Example.} Infinite monodromy.

Let $W, X_1, \ldots, Y_n, Z_0$ be as in \S 2 with $[X_i,Y_i] = Z_0$, and $\mathfrak{p} = span \{ X_1, \ldots, Y_n \}$. Let $a$ be the automorphism of $\mathfrak{h}$ fixing $Z_0$ and acting on each $span \{ X_i, Y_i \}$ by
$$
\left[
\begin{array}{lr}
 2 &  1 \\
1  &  1
\end{array}
\right]
$$

Let $b$ fix $Z_0$ and act on each $span \{ X_i, Y_i \}$ by 
$$
\left[
\begin{array}{lr}
 1 &  1 \\
\frac{-1 + \sqrt{5}}{2}  &  \frac{-1 - \sqrt{5}}{2}
\end{array}
\right]
$$

The conjugate $b^{-1} a b$ is diagonal.  Let $V$ fix $Z_0$ and act on each $span \{ X_i, Y_i \}$ by
$$
b \cdot
\left[
\begin{array}{lr}
 \ln( \frac{3 + \sqrt{5}}{2} ) &  0 \\
0  &  \ln( \frac{3 - \sqrt{5}}{2})
\end{array}
\right]
\cdot b^{-1}
$$

so $e^V = a$.  Define a path of automorphisms $\varphi_t = e^{tV}$, so $\nu_t = \left. V \right|_{\mathfrak{p}}$ for all $t$.

Obviously $\nu_t$ is infinitesimally symplectic; one may compute it is infinitesimally definite, so $\varphi_t$ is metric-defining.  Let $\Phi_t$ be the path of automorphisms of $H$ with derivative $\varphi_t$.  Since $\Phi$ is actually a homomorphism from $\R$ to $\mbox{Aut}(H)$, it is in particular $\Z$-equivariant.  Let $A = \Phi_1$, the automorphism with derivative $a$.  Let $\Lambda$ be the $\Z$-span of $X_1, \ldots, Y_n, \frac{1}{2}Z_0$ in $\mathfrak{h}$.  The Campbell-Hausdorff formula (see \cite{Helg} II.1.8) for an order-two nilpotent group reads
$$ e^X \cdot e^Y = e^{X+Y+ \frac{1}{2}[X,Y]}$$

The exponential of $\Lambda$ is a lattice subgroup $L$ in $H$.  Now $a^n(\Lambda) \subset \Lambda$ implies $A^n(L) \subseteq L$ for all $n \in \Z$.  The group $\Z \ltimes_A L$ acts isometrically and cocompactly on $\R *_{\Phi} H$ endowed with the Lorentz metric constructed above.  The quotient manifold $M$ admits an isometric $H$-action.  The quotient of $M$ by this action fibers over $S^1$, and the monodromy of this fiber bundle is the infinite group $<a>$.  Any quotient of a warped Heisenberg group $S$ by a lattice $\Gamma$ has finite monodromy, so $M$ has no finite cover diffeomorphic to such a quotient. 

Results of \S 4.3 say that all codimension-one $H$-actions arise from the construction in this section, assuming the actions are real-analytic.

\section{Toward classification of $H$-actions}

   We consider real-analytic actions of arbitrary codimension in Section 4.1; we prove a refinement of a result of Gromov to show that if the action of the center factors through a compact quotient, then almost every orbit map is proper for the lifted $H$-action on the universal cover.  In the second section, techniques of Gromov and Zimmer are combined with Lemma \ref{lemma} to show in the codimension-one case that the action on the universal cover is proper and free.  Theorem \ref{difftype} determines the diffeomorphism type of compact Lorentz manifolds admitting codimension-one, real-analytic, isometric $H_n$-actions.  Finally, we show that all such actions arise from the construction in Section 3.3; a complete classification follows.

The following two important properties of isometric actions of nilpotent groups are proved by Adams and Stuck in \cite{AS}.  They will be used extensively below.

\begin{thm}[Adams and Stuck \cite{AS} 6.8] \label{ASlocfree} 
If a Heisenberg group $H$ acts locally faithfully on a compact Lorentz manifold by isometries, then the $H$-action is locally free.
\end{thm}   

\begin{thm}[Adams and Stuck \cite{AS} 7.4] \label{ASlightlike}
For $H$ as above, the pullback to $\mathfrak{h}$ of the metric from any tangent space is Ad $H$-invariant and lightlike with kernel $\mathfrak{z}(\mathfrak{h})$.
\end{thm}

Recall that the action of a group $G$ on a manifold $M$ is \emph{proper} if for all compact subsets $A \subseteq M$, the set $G_A := \{ g \in G \ | \ gA \cap A \neq \emptyset \}$ is compact.  A group acting properly on a compact manifold is compact.  Orbits for a proper action are closed; moreover, orbits of any closed subgroup of a group acting properly are closed.  A map $f: X \rightarrow Y$ is \emph{proper} if $f^{-1}(K)$ is compact for all compact subsets $K$ of $Y$.

\subsection{Higher codimension: if $Z(H)$ factors through $S^1$}

Suppose the action of the center $Z(H)$ on $M$ factors through a proper action.  In \cite{Ze2} 4.3, Zeghib asks whether this factoring occurs for all isometric actions of Heisenberg groups on compact Lorentz manifolds.  Results of the next section imply that it does in the codimension-one case.

Let $\Gamma = \pi_1(M)$.   Let $\overline{\Lambda}$ be the kernel of the $H$-action on $M$.  This discrete normal subgroup must be central and, by assumption, cocompact in $Z(H)$.  The group of lifts of the $H$-action to $\widetilde{M}$ is an extension of $H / \overline{\Lambda}$ by $\Gamma$ with identity component a covering group of $H / {\overline{\Lambda}}$.  In other words, there is some subgroup $\Lambda \subseteq \overline{\Lambda}$ such that the identity component of the group of lifts is $H / \Lambda$, and $(H / \Lambda) \cap \Gamma \cong \overline{\Lambda}/ \Lambda$. 

At this point, it is necessary to assume that the Lorentz metric on $M$ is real-analytic, in order to apply the following theorem of Gromov.  For the case in which $H$ is simple, Zimmer gives a more detailed proof in \cite{Zi2} Corollary 4.3.  

\begin{thm}[\cite{Gromov} 5.2.A2] 
\label{ccontainsh}
Let $M$ be a compact real-analytic Lorentz manifold with an analytic locally free action of a connected Lie group $H$ by isometries.  Assume that the Zariski closure of Ad $H$ in Aut $\mathfrak{h}$ has no proper normal, algebraic, cocompact subgroup.  Let $\mathfrak{c}$ be the Lie algebra of global Killing fields on the universal cover $\widetilde{M}$ commuting with $\mathfrak{h}$.  For almost every $x \in \widetilde{M}$ and every $Y \in \mathfrak{h}$, there is $C \in \mathfrak{c}$ for which $C(x) = Y(x)$.
\end{thm}

This theorem in fact applies whenever $H$ acts locally freely and real-analytically preserving a finite measure and a real-analytic rigid geometric structure of algebraic type on $M$.  A real-analytic isometric Heisenberg action satisfies the hypotheses of the theorem.  

Next we show that almost every orbit on $\widetilde{M}$ is proper.   This result is an extension of a theorem of Gromov (\cite{Gromov} 6.1.A) to the case of $Z(H)$ non-compact, and is seemingly alluded to in \cite{Gromov} 6.1.C(b).  Gromov's result appears with proof for the case of $H$ simple in Zimmer's paper \cite{Zi2}.  

\begin{thm}
\label{hproper}
Let $M$ be a compact real-analytic Lorentz manifold with an analytic locally faithful action of the Heisenberg group $H$ such that the kernel of the action contains a cocompact subgroup of $Z(H)$.  Let $H / \Lambda$ be the identity component of the group of lifts to $\widetilde{M}$.  Then for almost every $x \in \widetilde{M}$, the orbit map $H / \Lambda \rightarrow (H/ \Lambda) . x \subseteq \widetilde{M}$ is proper.  
\end{thm}

\begin{Pf}
By Theorem \ref{ASlocfree}, the $H$-action on $M$ is locally free.  Let $T^{\mathcal{O}}(M)$ and $T^{\mathcal{O}}(\widetilde{M})$ be the subbundles of $T(M)$ and $T(\widetilde{M})$, respectively, tangent to $H$-orbits.  Let $\mathcal{F}^{\mathcal{O}}(M)$ and $\mathcal{F}^{\mathcal{O}}(\widetilde{M})$ be the corresponding bundles of frames.  These bundles are smoothly trivializable using the maps 
$$f_x(X) = \left. \frac{\partial}{\partial t} \right|_0 e^{tX}.x$$

for $X \in \mathfrak{h}$.  These maps satisfy the equivariance relation
$$f_{gx} \circ \mbox{Ad} g^{-1} = g_{*x} \circ f_x$$

Therefore, the $H$-action in either trivialization is $h(x, Y) = (hx, \mbox{Ad} h^{-1} (Y))$.  Since Ad : $H / \overline{\Lambda} \rightarrow \mbox{Ad}(H)$ is proper, $H/ \overline{\Lambda}$ is proper on $\mathcal{F}^{\mathcal{O}}(M)$.  

Let $\mathfrak{c}$ be as in Theorem \ref{ccontainsh} above.  The group of isometries of $\widetilde{M}$ leaving $\mathfrak{c}$ invariant acts on $\widetilde{M} \times \mathfrak{c}$ by $g. (x, C) = (g . x, g_* C)$, such that the evaluation map $\widetilde{M} \times \mathfrak{c} \rightarrow T (\widetilde{M})$ is equivariant.  In particular, $H/ \Lambda$ and $\Gamma$ have commuting actions on $\widetilde{M} \times \mathfrak{c}$, and these actions agree on the intersection $\overline{\Lambda}/ \Lambda$.  Let $q : \widetilde{M} \times \mathfrak{c} \rightarrow T(M)$ be the evaluation map composed with the quotient $T(\widetilde{M}) \rightarrow T(M)$.  Note that $q$ is equivariant with respect to $H / \Lambda \rightarrow H / \overline{\Lambda}$.

Suppose $x$ satisfies the conclusion of Theorem \ref{ccontainsh}, and suppose that $h .x \in K$ for some $h$ in $H/ \Lambda$ and some compact set $K$ in $\widetilde{M}$.  The goal is to find a compact subset $A \subset H / \Lambda$ containing all such $h$.  For this purpose, we may assume $(H/ \Lambda).x \cap K$ is dense in $K$.  

For the chosen $x$, there are $C_1, \ldots, C_m \in \mathfrak{c}$ such that $C_1(x), \ldots, C_m(x)$ are a framing of $T^{\mathcal{O}}_x(\widetilde{M})$.  Since each $C_i$ is $(H / \Lambda)$-invariant, $C_1(y), \ldots, C_m(y)$ are a framing of $T^{\mathcal{O}}_y(\widetilde{M})$ for all $y$ in the orbit of $x$.  Now suppose $y = \lim(h_n. x)$ for a sequence $h_n \in H / \Lambda$.  Each $C_i$ is smooth, and $T^{\mathcal{O}}(\widetilde{M})$ is a closed subbundle of $T \widetilde{M}$, so $C_i(y)$ is also in $T^{\mathcal{O}}(\widetilde{M})$ for $i = 1, \ldots, m$.  Suppose $C_i(y) = {\bf 0}$ for some $i$.  Now $C_i(y) = \lim C_i(h_n .x) = \lim h_{n*} C_i(x)$.  Let $C \in \mathfrak{h}$ be the image of $C_i(x)$ in the trivialization of $T^{\mathcal{O}}(\widetilde{M})$; then $\lim \mbox{Ad}(h_n^{-1}) C = {\bf 0}$.  Since $\mbox{Ad}(H) \subset SL(\mathfrak{h})$, there is some $C^{\prime} \in \mathfrak{h}$ with $\mbox{Ad}(h_n^{-1}) C^{\prime} \rightarrow \infty.$  But this $C^{\prime}$ also corresponds to some vector field in $\mathfrak{c}$, which is bounded on $K$, a contradiction.  Therefore, each $C_i$ is nonzero at all points of $K$, and $C_1(y), \ldots, C_m(y)$ are a framing of $T^{\mathcal{O}}_y(\widetilde{M})$ for all $y \in K$.

For each $y \in K$, the images $q(y, C_1), \ldots, q(y, C_m)$ form a frame on $T_y^{\mathcal{O}}(M)$.  The images $K \times \{ C_1 \}, \ldots, K \times \{ C_m \} $ form a compact subset $B$ of $\mathcal{F}^{\mathcal{O}}(M)$.  Since  $h (x, C_i) \in K \times \{ C_i \}$ for all $i$, the image of $h$ in $H / \overline{\Lambda}$ carries the frame $(q(x, C_1), \ldots, q(x, C_m))$ into $B$.  Since $H / \overline{\Lambda}$ is proper on $\mathcal{F}^{\mathcal{O}}(M)$, the images of all such $h$ are contained in a compact subset $\overline{A}$ of $H / \overline{\Lambda}$.

Let $A$ be a compact subset of $H / \Lambda$ mapping onto $\overline{A}$, and write $h = a \lambda $ with $\lambda \in \overline{\Lambda} / \Lambda$ and $a \in A$.  Now $h .x \in K$ implies $\lambda .x \in A^{-1} K$.  Since $\overline{\Lambda} / \Lambda$ is proper on $\widetilde{M}$, the set of all possible $\lambda$ is finite.  Thus the subset of $H / \Lambda$ carrying $x$ into a fixed compact set is compact, and the orbit map for $x$ is proper. 
\end{Pf}

{\bf Remark}: The conclusion of Theorem \ref{hproper} actually holds for any unimodular $H$ to which \ref{ccontainsh} applies, provided the $Z(H)$-action factors through a torus.

\subsection{Codimension-one: complete classification}

Throughout this section, $M$ is a compact real-analytic Lorentz manifold with a real-analytic locally faithful action of the Heisenberg group $H= H_n$ such that $\mbox{dim} M =2n+2$.  In this case, Gromov's centralizer theorem plus the basic fact of Lemma \ref{lemma} can be used to show properness of $H/ \Lambda$ on $\widetilde{M}$.   Two main results are that the quotient of $M$ by the $H$-action is a fiber bundle over $S^1$, and $\widetilde{M}$ is $\pi_1(M)$-equivariantly diffeomorphic to $\R \times H$.  We conclude by showing that all of these actions arise from the construction given in Section 3.3.

\begin{thm} 
\label{codim1hproper}
 Let $H / \Lambda$ be the identity component of the $H$-action lifted to $\widetilde{M}$.  Then $H/ \Lambda$ preserves a framing of $T(\widetilde{M})$; in particular, the $H / \Lambda$-action is free and proper.
\end{thm}

\begin{Pf}
  Let $\mathfrak{c}$ be the centralizer of $\mathfrak{h}$ on $\widetilde{M}$, as in the previous section, and let $x$ satisfy the conclusion of Theorem \ref{ccontainsh}.  Any $C$ in the kernel of the evaluation map $\mathfrak{c} \rightarrow T_x(\widetilde{M})$ integrates to a (local) isometry fixing $x$ with derivative fixing $f_x(\mathfrak{h})$ pointwise.  Such an isometry fixes a point and a frame by Lemma \ref{lemma}, and so is trivial (see \cite{Ko} II.1.3).  Thus the evaluation map on $\mathfrak{c}$ is injective, and by Theorem \ref{ccontainsh}, the dimension of $\mathfrak{c}$ is either $2n+2$ or $2n+1$.  

In the first case, $\mathfrak{c}(x) = T_x(\widetilde{M})$ for all $x \in \widetilde{M}$, and the evaluation map $\widetilde{M} \times \mathfrak{c} \rightarrow T \widetilde{M}$ defines an $H / \Lambda$-equivariant map $\widetilde{M} \times \mathcal{F}(\mathfrak{c}) \rightarrow \mathcal{F}(\widetilde{M})$, where $\widetilde{M} \times \mathcal{F}(\mathfrak{c})$ is the bundle of frames on $\widetilde{M} \times \mathfrak{c}$. The trivial $H / \Lambda$-map $\widetilde{M} \rightarrow \widetilde{M} \times \{ \omega \}$, for any $\omega \in \mathcal{F}(\mathfrak{c})$, gives an $H / \Lambda$-equivariant section $\widetilde{M} \rightarrow \mathcal{F}(\widetilde{M})$. 

In the second case, the set of $x$ such that $f_x(\mathfrak{h}) = \mathfrak{c}(x)$ is closed and has full measure, so it must be all of $\widetilde{M}$.  Fix a generator $Z_0$ for $\mathfrak{z}(\mathfrak{h})$, and note that $Z_0 \in \mathfrak{c}$.  Let $\widetilde{M} \times \mathcal{F}(\mathfrak{c})_0 \subset \widetilde{M} \times \mathcal{F}(\mathfrak{c})$ be the subbundle of frames with first element the Killing field for $Z_0$.  In any such frame, the remaining vectors evaluated at any $x \in \widetilde{M}$ span a spacelike subspace $\mathfrak{p}$ by Theorem \ref{ASlightlike}.  Take $W$ to be the unique transverse lightlike vector field that everywhere has inner product one with $Z_0^*$ and is orthogonal to $\mathfrak{p}$.  An element of $\mathcal{F}(\mathfrak{c})_0$ thus determines a frame on $T_x(\widetilde{M})$, and there is a well-defined map $\widetilde{M} \times \mathcal{F}(\mathfrak{c})_0 \rightarrow \mathcal{F}(\widetilde{M})$ that is smooth and $H / \Lambda$-equivariant. As above, composition with $\widetilde{M} \rightarrow \widetilde{M} \times \{ \omega \}$ for any $\omega \in \mathcal{F}(\mathfrak{c})_0$ yields an $H / \Lambda$-equivariant section. 
\end{Pf}

\begin{thm}
\label{difftype}
 There is a diffeomorphism of the universal cover $\widetilde{M}$ with $\R \times H_n$ carrying $H_n$-orbits to $H_n$-factors.  There is a lattice $L < H_n$ and a homomorphism $\Psi : \Z \rightarrow \mbox{Aut}(H)$ such that $\pi_1(M) \cong \Z \ltimes_{\Psi} L$. The lattice $L$ acts along the $H_n$-factors of $\widetilde{M}$, and $\Z$ acts by translations along $\R$.
\end{thm}

\begin{cor}
The quotient $H_n \backslash M$ is a fiber bundle over $S^1$, giving rise to a real-analytic map $S^1 \rightarrow PD(\mathfrak{h} / \mathfrak{z})$, the space of positive-definite inner products on $\mathfrak{h} / \mathfrak{z}$.
\end{cor}

\begin{Pf} of corollary

If $H \backslash \widetilde{M} \cong \R$, then $H \backslash M \cong H \backslash \widetilde{M} / \Gamma \cong \R / \Z \cong S^1$. 

By Theorem \ref{ASlightlike}, the pullback of the metric from each $H$-orbit to $\mathfrak{h}$ is lightlike and descends to a positive-definite inner product on $\mathfrak{h} / \mathfrak{z}$.
\end{Pf}

\begin{Pf} of theorem

We will denote by $\mathfrak{c}$ the Lie algebra of the centralizer of $\mathfrak{h}$ on $\widetilde{M}$ and by $\Gamma$ the fundamental group $\pi_1(M)$.

{\bf step 1}: $H$ acts faithfully.

Let $\Lambda$ be the kernel of the action on $\widetilde{M}$, as in the previous theorem.  Since $H / \Lambda$ acts properly and freely by Theorem \ref{ASlocfree}, $\widetilde{M}$ is a principal bundle over a connected $1$-manifold.  The fibers are connected, so the quotient is simply connected, hence, diffeomorphic to $\R$.  From the long exact sequence of homotopy groups for this fiber bundle, one computes that the fibers must be simply connected.  Thus $\Lambda$ is trivial, and $H$ acts faithfully. 

{\bf step 2}: $C$ acts properly and freely.

 Pick $Z_0 \in \mathfrak{z}(\mathfrak{h})$ and the rest of a basis $Z_0, Y_1, \ldots, Y_{2n}$ for $\mathfrak{h}$.  Let $\mathfrak{p}= span \{ Y_1, \ldots, Y_{2n} \}.$  Let $W$ be the transverse vector field defined with respect to $Z_0$ and $\mathfrak{p}$ by ({\bf M1}) - ({\bf M3}). The Killing fields $Z_0^*, Y_1^*, \ldots, Y_{2n}^*$ together with $W$ form a framing on $\widetilde{M}$ that descends to $M$.  This framing of $\widetilde{M}$ is also invariant by the flow along any Killing field in $\mathfrak{c}$.  It follows that each Killing field in $\mathfrak{c}$ is complete.  Let $C$ be the centralizer in Isom($\widetilde{M}$) of $H$.  The Lie algebra of $C$ is $\mathfrak{c}$, and $\Gamma$ lies in $C$.  Because it preserves a framing on $\widetilde{M}$, the group $C$ acts properly and freely.  

{\bf step 3}: if $\mbox{dim} \mathfrak{c} = 2n+2$, then $\mathfrak{c} \cong \mathfrak{s}$.

Suppose dim $\mathfrak{c} = 2n+2$, so $C^0$-orbits are open.  They are closed because $C$ acts properly; therefore, $C$ acts transitively.  Since the action is free, $C=C^0$, and $\Gamma$ is a lattice in $C$.

 Pick any $x_0$, and let $G$ be the set of elements $g$ of $C$ such that $g x_0 = h x_0$ for some $h \in H$.  Since $H$ acts freely on $\widetilde{M}$, there is a well-defined map $\theta_0 : G \rightarrow H$ satisfying $\theta_0(g) g. x_0 = x_0$.  This map is injective because $C$ acts freely, and it respects multiplication because $H$ and $C$ commute.  Therefore $\theta_0$ is an isomorphism from $G$ to $H$. 

The orbit map $C \rightarrow C. x_0 \subset \widetilde{M}$ intertwines the $G \times G$-action on $C$ by left and right translation, respectively, with the $G \times G$-action on $\widetilde{M}$ by the given $G$-action and $\theta_0 \circ \iota$, respectively, where $\iota$ is inversion in $G$.
 
The orbit map also yields an identification of $T_{x_0} \widetilde{M}$ with $\mathfrak{c}$.  Denote by $<,>$ the pullback of the Lorentz metric at $x_0$ to $\mathfrak{c}$.  The adjoint represention of $G$ on $\mathfrak{c}$ is equivalent to the $G$-action by $\theta(g)_* g_*$ on $T_{x_0} \widetilde{M}$.  Thus $<,>$ is infinitesimally $\mathfrak{g}$-invariant---that is, for any $Y \in \mathfrak{g}$ and $X, W \in \mathfrak{c}$, 
$$<[Y,X],W> + <X, [Y,W] > = 0$$  
Let $W$ be a vector complementary to $\mathfrak{g}$ in $\mathfrak{c}$.  For any $Y, Y^{\prime}$ in $\mathfrak{g}$, 
\begin{eqnarray*}
<[W,Y], Y^{\prime}>  +  <Y, [W,Y^{\prime}] >   & =  & - <[Y,W], Y^{\prime} >  -   < Y, [Y^{\prime},W]>  \\
  &  =  & <W, [Y,Y^{\prime}]>  +  <[Y^{\prime},Y],W>  \\
  &  =  &   0
\end{eqnarray*} 
For any $Y \in \mathfrak{g}$, 
\begin{eqnarray*}
 <[W,Y],W> & = & - <[Y,W],W> \\
                        & = & <W,[Y,W]> \\
                        & = & 0 
\end{eqnarray*} 
Since $\mathfrak{g}$ and $W$ together span $\mathfrak{c}$, it follows that $<,>$ is infinitesimally $W$-invariant.  Then $<,>$ is $\mbox{Ad}C$-invariant, and there is an isometric $C \times C$-action on $\widetilde{M}$ extending the $C$-action.  This action descends to a transitive isometric $C$-action on $M$.  By the classification of compact homogeneous Lorentz manifolds \cite{Ze1}, $\widetilde{M}$ is equivariantly isometric to $\widetilde{S}_{\lambda}$ for some $\lambda \in \Q^n$, and $\Gamma \cong \Z \ltimes_{\Psi} L$ for $L$ a lattice in $H_n$. The homomorphism $\Psi(q) = Ad(e^{qW})$ and has finite cyclic image.

Assume now that $\mbox{dim} \mathfrak{c} = 2n+1$.

{\bf step 4} : lattice in $H$.

Let $C^{\prime}$ be the closed subgroup generated by $C^0$ and $\Gamma$.  Properness of $C$ on $\widetilde{M}$ from step {\bf 2} implies that $C^0$-orbits have closed image in $\widetilde{M} / \Gamma = M$. These are diffeomorphic to $C^0/(C^0 \cap \Gamma)$, so $L = C^0 \cap \Gamma$ is a lattice in $C^0$.  The set of $x$ on which $f_x(\mathfrak{h})$ and $\mathfrak{c}(x)$ coincide is closed and full measure, so it is all of $\widetilde{M}$.  Then $C^0$-orbits coincide with $H$-orbits.  As in step {\bf 3}, choosing $x_0$ determines an isomorphism $\theta_0 : C^0 \rightarrow H$.  Then $\theta_0 (L)$ is a lattice in $H$.

{\bf step 5}: tranvserse foliation.

Let $Z_0, Y_1, \ldots, Y_{2n}$ be the basis for $\mathfrak{h}$ from step {\bf 2}; let $\mathfrak{p}$ be the span of $Y_1, \ldots, Y_{2n}$.  Let $W$ be the vector field on $\widetilde{M}$ satisfying ({\bf M1}) - ({\bf M3}) with respect to $\mathfrak{p}$ and $Z_0$; it is a $C$-invariant vector field transverse to $H$-orbits.  Since $W$ descends to $M$, it is complete.  Let $x_0$ be as in step {\bf 4}, and let $\mathcal{W}$ be an integral curve for $W$ with $\mathcal{W}(0) = x_0$.

Let $m$ be a right-invariant Riemannian metric on $C^0$.  Define a Riemannian metric on each $C^0$-orbit by identifying Killing fields with right-invariant vector fields on $C^0$.  Now define a Riemannian metric on $\widetilde{M}$ by making $W$ norm 1 and orthogonal to $C^0$-orbits.  The product metric on $\R \times C^0$ arising from $\frac{\partial}{\partial t}$ and $m$ is complete.  The map $(t,g) \mapsto g.\mathcal{W}(t)$ is a local isometry, so it is a Riemannian covering $\R \times C^0 \rightarrow \widetilde{M}$ (see \cite{BH}, \cite{CG}, \cite{Gromov} {\bf 5.3.D}, \cite{Ze1} \S 5.5 for origins and more ambitious applications of this idea).  Also, $C^0$-orbits equal $H$-orbits, so applying $\theta_0 \circ \iota$, where $\iota$ is inversion, on $C^0$ gives a diffeomorphism of $\widetilde{M}$ with $\R \times H$ carrying $H$-orbits to factors $\{ t \} \times H$.  Since $C$ preserves both the vector field $W$ and the foliation by $H$-orbits, the $C$-action, and, in particular, the $\Gamma$-action, respect the product structure induced on $\widetilde{M}$.  Denote by $t(x)$ and $p(x)$ the resulting coordinate functions to $\R$ and $H$, respectively.

{\bf step 6}: splitting of $C$.

Let $Q \subset C$ be the subgroup preserving $\mathcal{W}$.  Note that $c \in Q \Leftrightarrow c.x_0 \in \mathcal{W}$. In fact, because any $c \in C$ commutes with flow along $W$, any $c \in Q$ acts on $\mathcal{W}$ by translation by $t(c.x_0)$.  The subgroup $Q$ contains the translation by $t(c.x_0)$ for all $c \in C$ because there is always some $g \in C^0$ with $gc.x_0 \in \mathcal{W}$.

Since $C$ acts freely, it is clear that each element can be written as a product in $Q \cdot C^0$.  Define
\begin{eqnarray*}
\Psi & : &  Q \rightarrow \mbox{Aut}(C^0) \\
\Psi(q)(g) & = & q^{-1} g q
\end{eqnarray*}

 Note that $\Psi$ is conjugate by $\theta_0$ into $\mbox{Aut}(H)$.  We have $C \cong Q \ltimes_{\Psi} C^0$.

{\bf step 7}: splitting of $\Gamma$.

 The quotient $C^0 \backslash \widetilde{M}$ can be identified with $\mathcal{W} \cong \R$, and $\Gamma / L$ acts faithfully on it with quotient $C^{\prime} \backslash \widetilde{M} \cong S^1$.  Thus the projection of $\Gamma$ to $Q$ is the group of translations by $\alpha \Z$ for some $\alpha \neq 0$. By the splitting of $C$ from the previous step, $\Gamma$ is isomorphic to a semi-direct product $\Z \ltimes_{\Psi} L$.
\end{Pf}

{\bf Remark 1:} Let $\theta_0 : C^0 \rightarrow H$ be as in step {\bf 5}.  Under the diffeomorphism $\widetilde{M} \cong \R \times H$, the action of $q \in \Z$ is
$$q.x = (t(x) + q \alpha,(\theta_0 \circ \Psi_q \circ \theta_0^{-1})(p(x))$$

{\bf Remark 2:} The group $C$ inherits a real-analytic structure from $\widetilde{M}$ for which it acts real-analytically.  The isomorphism $\theta_0$ is real-analytic.  The vector field $W$ is real-analytic, as are the curve $\mathcal{W}$ and the diffeomorphism $\widetilde{M} \rightarrow \R \times H$.

In the remaining statements, we will say, given a metric-defining path $\Phi_t$, that a space $X$ is equivariantly isometric to $\R *_{\Phi} H$ if there is some $Z_0 \in \mathfrak{z}$ for which $X$ is equivariantly isometric to $\R *_{\Phi} H$ in the metric determined by $\Phi_t$ and $Z_0$. 

\begin{thm}
\label{metrictype}
The universal cover $\widetilde{M}$ is $H_n$-equivariantly isometric to $\R *_{\Phi} H_n$ for a real-analytic,  metric-defining, lattice compatible path $\Phi : \R \rightarrow \mbox{Aut}(H_n)$. The restriction of $\Phi$ to $\Z$ coincides with the homomorphism $\Psi$ of the previous theorem.
\end{thm}

\begin{Pf}  

Let $\mathfrak{c}$ and $C$ be as in step {\bf 2} of the previous proof.

{\bf step 1}: homogeneous case.

If $\mbox{dim} \mathfrak{c} = 2n+2$, then $M$ is homogeneous by step {\bf 3} in the previous proof.  By the results of \cite{Ze1}, there is $\lambda \in \Q^n$ for which $C$ is isometric to $\widetilde{S}_{\lambda}$ endowed with the unique bi-invariant metric, and $M$ is equivariantly isometric to $C / \Gamma$ for $\Gamma$ a lattice subgroup. The group $\widetilde{S}_{\lambda}$ is $H$-equivariantly isometric to $\R *_{\Phi} H$, where $\R$ parametrizes $e^{tW}$ and $\Phi_t(h) = e^{-tW}he^{tW}$.  The path $\Phi_t$ is a homomorphism, so equivariance is clear.  Note $\Phi$ restricted to $\Z$ equals $\Psi$ from step {\bf 3} of the previous theorem.

Now assume $\mbox{dim} \mathfrak{c} = 2n+1$.

{\bf step 2}: metric-defining path.

Let $\mathfrak{h} = \mathfrak{p} \oplus \R Z_0$ be the decomposition in the proof of Theorem \ref{difftype} and $W$ the transverse vector field with integral curve $\mathcal{W}$ corresponding to this decomposition as above.

Let $\theta_t : C^0 \rightarrow H$ be the isomorphism defined by $\theta_t(g)g.\mathcal{W}(t) = \mathcal{W}(t)$; these form a real-analytic path of real-analytic isomorphisms.  Let $\Phi_t = \theta_0 \circ \theta_t^{-1}$; these form a real-analytic path.  Denote by $(t,p)$ the coordinates of the diffeomorphism $\widetilde{M} \rightarrow \R \times H$. The $C^0$-action in these coordinates is $g.x = (t(x),\theta_0(g^{-1})\cdot p(x))$.  It follows that $h.x = (t(x), \Phi_{t(x)}(h) \cdot p(x))$, so $\widetilde{M}$ is $H$-equivariantly diffeomorphic to $\R *_{\Phi} H$.

Fix a point $x$ and let $t(x) = t$ and $p(x) = p$.  The vector field $W$ is carried to $\frac{\partial}{\partial t}$ on $\R *_{\Phi} H$, and Killing fields are preserved.  Both $\frac{\partial}{\partial t}$ and $W$ have properties ({\bf M1})-({\bf M3}) with respect to $\mathfrak{p}$ and $Z_0$.  Also, $Z_0^*$ is isotropic and orthogonal to $f_x(\mathfrak{p})$ in both metrics.  It remains to show that $\Phi_t$ is metric-defining and that our map is isometric restricted to $f_x(\mathfrak{p})$.  As in the proof of Proposition \ref{ConstrIsIsometric}, for $X \in \mathfrak{h}$, and $X^*$ the corresponding Killing field,
$$[W,X^*](x) = ( \varphi_t^{-1} \circ \left( \frac{\partial \varphi_s}{\partial s} \right)_t(X))^*(x)$$
 
where $\varphi_t = D_e(\Phi_t)$.
Since $Z_0^* \in \mathfrak{c}$, 
$$[W,Z_0^*](x) = {\bf 0}$$

Also as in Proposition \ref{ConstrIsIsometric}, the fact that $H$ acts isometrically implies that for all $K \in \mathfrak{h}$ and $x \in \widetilde{M}$,
$$<[K^*,W],Y>_x + <W,[K^*,Y]>_x = 0$$

where $Y$ is any vector field such that $<W,Y>$ is constant along $H$-orbits.  Taking $K \in \mathfrak{p}$ and $Y = W$ gives
$$2<[K^*,W],W>_x = 0$$

Since $[K^*,W](x)$ belongs to $f_x(\mathfrak{h})$, it is also orthogonal to $Z_0^*(x)$.  Thus $[K^*,W](x)$ always belongs to $f_x(\mathfrak{p})$.

It follows that $\varphi_t^{-1} \circ \left( \frac{\partial \varphi_s}{\partial s} \right)_t $ annihilates $Z_0$ and preserves $\mathfrak{p}$ for all $t$.
Let $\nu_t$ be the restriction of this endomorphism to $\mathfrak{p}$. 

Take $Y = Y^*$ and $K^*$ both Killing fields from $\mathfrak{p}$.  Recall that $X \mapsto X^*$ is a Lie algebra homomorphism.  Then 
$$<[K^*,W],Y^*>_x + <W,[K^*,Y^*]>_x = <[K^*,W],Y^*>_x + \omega_0(K,Y)$$

because $[K,Y] = \omega_0(K,Y)$ and $<W,Z_0^*> = 1$ everywhere.  Thus
$$<[W,K^*],Y^*>_x = <\nu_t(K)^*,Y^*>_x = \omega_0(K,Y)$$

Because $\omega_0$ is non-degenerate, the kernel of $\nu_t$ is trivial.  We obtain
$$<K^*,Y^*>_x = \omega_0(\nu_t^{-1}K,Y)$$

Now symmetry of the metric implies $\nu_t^{-1}$ and $\nu_t$ are infinitesimally symplectic, and definiteness implies $\nu_t$ is infinitesimally definite for all $t$.  Thus $\Phi_t$ is metric-defining, and the map $\widetilde{M} \rightarrow \R *_{\Phi} H$ is isometric.

{\bf step 3}: equivariance of $\Phi$

Recall $\Gamma = \pi_1(M) \subset C$ and the splitting $\Gamma \cong \Z \ltimes_{\Psi} L$ from step {\bf 6} of the previous proof.
We relate $\Phi_{q \alpha}$ and $\Psi_q$ for $q$ in the projection of $\Gamma$ to $\Z$.  Applying remark {\bf 1} above,
$$(q \circ h).x_0 = (q \alpha,(\theta_0 \circ \Psi_q \circ \theta_0^{-1})(h))$$

while
$$(h \circ q). x_0 = (q \alpha,\Phi_{q \alpha}(h))$$

Because $H$ and $\Gamma$ commute, we conclude $\Phi_{q \alpha} = \theta_0 \circ \Psi_q \circ \theta_0^{-1}$.

Next, we show equivariance.  Again applying remark {\bf 1} for arbitrary $x$ with coordinates $t(x) = t$ and $p(x) = p$,
$$(h \circ q).x = (t + q \alpha,\Phi_{q \alpha +t}(h) \cdot \Phi_{q \alpha}(p))$$

This must equal
$$(q \circ h). x = (t + q \alpha, \Phi_{q \alpha}(\Phi_t(h) \cdot p)) = (t + q \alpha, (\Phi_{q \alpha} \circ \Phi_t)(h) \cdot \Phi_{q \alpha}(p))$$

Therefore, $\Phi_{q \alpha +t} = \Phi_{q \alpha} \circ \Phi_t$.  

Take a new metric-defining path $\Phi_t^{\prime} = \Phi_{\alpha t}$. As in Section 3.3.2, there is an $H$-equivariant isometry $\R *_{\Phi} H \rightarrow \R *_{\Phi^{\prime}} H$, and now $\Phi^{\prime}$ is lattice-compatible, with $\Phi^{\prime}_q = \theta_0 \circ \Psi_q \circ \theta_0^{-1}$ for all $q \in \Z$.
\end{Pf}

\begin{thm}(Classification).
\label{classfication}
There is a bijective correspondence between 

(a)  universal covers of compact Lorentz manifolds with real-analytic codimension-one isometric $H_n$-actions, up to equivariant isometry

(b) real-analytic, lattice-compatible, metric-defining paths $\Phi_t$ in $\mbox{Aut}(H_n)$, up to equivalence of metric-defining paths.
\end{thm}

\begin{Pf}
Let $M$ be a compact Lorentz manifold as in (a). By Theorem \ref{metrictype}, $\widetilde{M}$ is $H$-equivariantly isometric to some $\R *_{\Phi} H$.  If $\widetilde{M}$ is $H$-equivariantly isometric to $\R *_{\Phi^{\prime}} H$ for another path $\Phi^{\prime}$ (arising from a different choice of splitting $\mathfrak{h} = \mathfrak{p}^{\prime} \oplus \R Z_0^{\prime}$), then Proposition \ref{HEquiv} implies $\Phi$ and $\Phi^{\prime}$ are equivalent.

Given a path $\Phi_t$ as in (b), let $L$ be a lattice such that $\Phi(q)$ preserves $L$ for all $q \in \Z$, and let $\Gamma = \Z \ltimes_{\Phi} L$.  Then $(\R *_{\Phi} H) / \Gamma$ is a compact Lorentz manifold as in (a) by Propositions \ref{ConstrIsIsometric} and \ref{GammaIsIsometric}.  If $\Phi_t^{\prime}$ is an equivalent metric-defining lattice-compatible path, then $\R *_{\Phi} H$ and $\R *_{\Phi^{\prime}} H$ are $H$-equivariantly isometric from Section 3.3.2.
\end{Pf}

\end{document}